\documentclass[11pt, psfig]{article}
\usepackage{amssymb,amsmath,amscd}
\newtheorem{lem}{Lemma}[subsection]
\newtheorem{conj}{Conjecture}
\newtheorem{thm}{Theorem}[subsection]
\newtheorem{cor}{Corallary}[subsection]
\newtheorem{defn}{Definition}[subsection]

\newcommand{\f}[1]{\mathfrak{#1}}

\newcommand{\p}{\prime}
\newcommand{\mb}{\mathbb}
\newcommand{\commentout}[1]{}

\newcommand{\mc}{\mathcal}
\newcommand{\arr}[1]{\left( \begin{array}{clcr} #1 \end{array} \right)}

\newcommand{\diag}{{\ \rm diag}}
\newcommand{\xin}{Sp_{2n}(\mb R)}

\begin{document}
\title{Compositions of Theta Correspondences}
\author{Hongyu He \\
Department of Mathematics \& Statistics
\footnote{email:hongyu@math.lsu.edu } 
\footnote{AMS Subject Primary 22E45, 22E46} }
\date{}
\maketitle
\abstract{Theta correspondence $\theta$ over $\mb R$ is established by Howe in 
(~\cite{howe}). In ~\cite{unit}, we prove that $\theta$
preserves unitarity under certain restrictions, generalizing the result of Jian-
Shu Li (~\cite{li2}). The goal of this paper is to elucidate the idea of 
constructing unitary representation through the propagation of theta 
correspondences. 
We show that under a natural condition on the sizes of the
 related dual pairs which can be predicted by the orbit method (~\cite{pdk}, 
~\cite{vogan3}, ~\cite{pan}), 
one can compose theta correspondences to obtain unitary representations. We call 
this process quantum induction.
\section{Introduction}
An important problem in representation theory is the classification and 
construction of irreducible unitary representations. Let $G$ be a reductive 
group and $\Pi(G)$ be its admissible dual. For an algebraic semisimple group $G$, the admissible dual $\Pi(G)$ is known mostly due to the works of Harish-Chandra, R. Langlands, and Knapp-Zuckerman (see ~\cite{langlands}, ~\cite{kz}).
Let $\Pi_u(G)$ be the set of equivalence classes of 
irreducible unitary representations of $G$, often called the unitary dual of 
$G$. The unitary dual of general linear groups is classified by Vogan (see 
~\cite{vogan0}). The unitary dual of complex classical groups is classified by 
Barbasch (see ~\cite{barbasch}). Recently, Barbasch has classified all the 
spherical duals for split classical groups (see ~\cite{barbasch2}). The unitary 
duals $\Pi_u(O(p,q))$ and $\Pi_u(Sp_{2n}(\mb R))$ are not known in general. \\
\\
In ~\cite{howesmall}, Roger Howe constructs certain small unitary 
representations of the symplectic 
group using Mackey machine. Later, Jian-Shu Li generalizes Howe's construction 
of small unitary representations to all classical groups. 
In particular, Li defines a sesquilinear form $(,)_{\pi}$ that relates these 
constructions to the
theta correspondence (see ~\cite{howe79}, ~\cite{li2}). It then becomes 
clear to many people that some irreducible unitary representations 
can be constructed through the 
propagation of theta correspondences (See ~\cite{li97}, ~\cite{hl} and 
~\cite{pr1} and the references within them).
So far, constructions 
can only be carried out for "complete small orbits" (see ~\cite{li97}). 
The purpose of this paper is to make it work for nilpotent orbits in general, 
for real orthogonal groups and symplectic groups. 
 \\ 
\\
Consider 
the group $O(p,q)$ and $Sp_{2n}(\mb R)$. The theta correspondence with respect 
to $$O(p,q) \rightarrow Sp_{2n}(\mb R)$$ is formulated by Howe as a one-to-one 
correspondence
$$\theta(p,q;2n): \mc R(MO(p,q),\omega(p,q;2n)) \rightarrow
\mc R(MSp_{2n}(\mb R), \omega(p,q;2n))$$
where $MO(p,q)$ and $MSp_{2n}(\mb R)$ are some double coverings of $O(p,q)$ and 
$Sp_{2n}(\mb R)$ respectively and
$$\mc R(MO(p,q),\omega(p,q;2n)) \subseteq \Pi(MO(p,q)); \qquad
\mc R(MSp_{2n}(\mb R), \omega(p,q;2n)) \subseteq \Pi(MSp_{2n}(\mb R))$$ 
(see ~\cite{howe}). We denote the inverse of $\theta(p,q;2n)$ by 
$\theta(2n;p,q)$. For the sake of simplicity, we define
$$\theta(p,q;2n)(\pi)=0  $$
if $\pi \notin \mc R(MO(p,q), \omega(p,q;2n))$. We define $\theta(p,q;2n)(0)=0$ 
and $0$ can be regarded as the NULL
representation. \\
\\
For example, given an "increasing" string 
$$O(p_1,q_1) \rightarrow Sp_{2n_1}(\mb R) \rightarrow O(p_2,q_2) \rightarrow 
Sp_{2n_2}(\mb R) \rightarrow \ldots \rightarrow Sp_{2n_m}(\mb R) \rightarrow 
O(p_m,q_m) $$
$$p_1+q_1 \equiv p_2+q_2 \equiv \ldots \equiv p_m+q_m \qquad (\mod 2),$$
consider the propagation of theta correspondence along this string:
$$\theta(2n_m; p_m,q_m) \ldots \theta(2n_1;p_2,q_2) \theta(p_1,q_1;2n_1)(\pi).$$
Under some favorable conditions on $\pi \in \Pi_u(O(p_1,q_1))$, one hopes to 
obtain
a unitary representation in $\Pi_u(O(p_m,q_m))$. In this paper, we supply a 
sufficient condition for 
$$\theta(2n_m; p_m,q_m) \ldots \theta(2n_1;p_2,q_2) \theta(p_1,q_1;2n_1)(\pi)$$
to be unitary. We denote the resulting representation of $MO(p_m,q_m)$ by
$$Q(p_1,q_1; 2n_1; p_2,q_2;2n_2; \ldots ; p_m,q_m)(\pi).$$
We call $Q(p_1,q_1; 2n_1; p_2,q_2;2n_2; \ldots ; p_m,q_m)$ quantum induction.  
In addition to the assumption that certain Hermitian forms do not vanish, we 
must also assume the matrix coefficients of $\pi$ satisfy a mild growth 
condition. \\
\\
Based on the work of Przebinda (~\cite{pr}), we further determine the behavior 
of infinitesimal characters under quantum induction.
In certain limit cases, the infinitesimal character under quantum induction 
behaves exactly in the same way as under parabolic induction. In fact, in some 
limit cases, quantum induced representations can be obtained from unitarity-
preserving parabolic induction (see ~\cite{quan}). Finally, motivated by the 
works of Przebinda and his collaborators, we
make a precise conjecture regarding the associated variety of the quantum 
induced representations (Conjecture 2). \\
\\
There is one problem we did not address in this paper, namely, the nonvanishing 
of certain Hermitian forms $(,)_{\pi}$ with $\pi \in \Pi(Mp_{2n}(\mb R))$. In a 
forthcoming article (~\cite{quan}), we partially address this problem and 
construct a set of special unipotent representations in the sense of Vogan (see 
~\cite{vogan89}). I wish to
thank Prof. Shou-En Lu for her encouragements and the referee for some very helpful comments. 
\section{Main Results}

\subsection{Notations}
In this paper, unless stated otherwise, all representations are regarded as 
Harish-Chandra modules. This should cause no problems since most
representations in this paper will be admissible with respect to a reductive 
group. Thus
unitary representations in this paper would mean unitrizable Harish-Chandra 
modules.  "Matrix coefficients" of a
representation $\pi$ of a real reductive group $G$ will refer to the $K$-finite 
matrix coefficients with respect to a maximal compact subgroup $K$. A vector $v$ in
an admissible representation $\pi$ means that $v$ is in the Harish-Chandra 
module of $\pi$ which shall be evident within the context.\\
\\
Let $(G_1,G_2)$ be a reductive dual pair of type I (see ~\cite{howe} 
~\cite{li2}).
The dual pairs in this paper will be considered as ordered. For example, the 
pair $(O(p,q), Sp_{2n}(\mb R))$ is considered different from the pair 
$(Sp_{2n}(\mb R), O(p,q))$. Unless stated otherwise,
we
will, in general, assume that the size of $G_1(V_1)$ is less or equal to the 
size of $G_2(V_2)$, i.e.,
$\dim_D(V_1) \leq \dim_D(V_2)$. Let $(G_1,G_2)$ be a dual pair in the symplectic 
group $Sp$. Let $Mp$ be the unique double covering of $Sp$. Let $\{1,\epsilon\}$ 
be the preimage of the identity element in $Sp$.  For a subgroup $H$ of $Sp$, 
let
$MH$ be the preimage of $H$ under the double covering. Whenever we use the 
notation $MH$, $H$ is considered to be a subgroup of certain $Sp$ which shall be 
evident within the context. Let $\omega(MG_1, MG_2)$ be a Schr\"odinger model of 
the oscillator representation of $Mp$ equipped with a dual pair $(MG_1,MG_2)$. 
The Harish-Chandra module of $\omega(MG_1, MG_2)$ consists of polynomials 
multiplied by the Gaussian function. 
Since the pair $(G_1,G_2)$ is ordered, we use $\theta(MG_1, MG_2)$ to denote the
theta correspondence from $\mc R(MG_1, \omega(MG_1, MG_2))$ to
$\mc R(MG_2, \omega(MG_1, MG_2))$. 
We use $\bold n$ to denote the constant vector
$$(n,n, \ldots, n)$$
The dimension of $\bold n$ is determined within the context.
Finally, we say a vector 
$$x=(x_1,x_2, \ldots x_n) \prec 0$$ if
$$\sum_{j=1}^k x_j < 0 \qquad \forall k \geq 1$$
and $x \preceq 0$ if
$$\sum_{j=1}^k x_j \leq 0 \qquad \forall k \geq 1.$$
\\
In this paper,
the space of $m \times n$ matrices will be denoted by $M(m,n)$. The set of 
nonnegative integers will be denoted
by $\mb N$. For the group $O(p,q)$, we assume that $p \leq q$ unless stated 
otherwise. For a reductive group $G$, $\Pi(G)$, $\Pi_u(G)$ will be the 
admissible dual and the unitary dual respectively. \\
\\
We extend the definition of matrix coefficients to the NULL representation.
The matrix coefficients of the NULL representation is defined to be the zero 
function.
\subsection{Theta Correspondence in Semistable Range and Unitary 
Representations}
Let $\pi \in \Pi(MG_1)$.
Following ~\cite{li2}, for every $u,v \in \pi$ and $\phi, \psi \in \omega(MG_1, 
MG_2)$, we formally define
\begin{equation}~\label{avera}
(\phi \otimes  v, \psi \otimes u)_{\pi}=\int_{MG_1}(\omega(MG_1, MG_2) (\tilde 
g_1) \phi, \psi)(u, \pi(\tilde g_1) v) d \tilde g_1.
\end{equation}
Roughly speaking, if  the functions
$$(\omega(MG_1, MG_2) (\tilde g_1) \phi, \psi)(u, \pi(\tilde g_1) v) \qquad 
(\forall \phi, \psi \in 
\omega(MG_1, MG_2); \forall u, v \in \pi) $$
are in $L^{1}(MG_1)$ and $\pi(\epsilon)=-1$, $\pi$ is said to be in the
semistable range of $\theta(MG_1, MG_2)$ (see ~\cite{theta}). We denote the semistable range
of $\theta(MG_1, MG_2)$ by $\mc R_s(MG_1, MG_2)$.  \\
\\
{\bf Suppose from now on} that $\pi \in \mc R_s(MG_1, MG_2)$. In ~\cite{theta}, 
we showed that if $(,)_{\pi}$ does not vanish, then $(,)_{\pi}$ descends into a 
Hermitian form
on $\theta(MG_1, MG_1)(\pi)$. 
For $\pi \in \mc R_s(MG_1, MG_1)$, we {\bf define }
\begin{equation}
\theta_s(MG_1,MG_1)(\pi)=\left\{ \begin{array}{cc}
 \theta(MG_1,MG_2)(\pi) & \mbox{ if $(,)_{\pi} \neq 0$} \\
 0 & \mbox{ if $(,)_{\pi}=0$} \end{array} \right.
\end{equation}
$\theta_s(MG_1, MG_2)(\pi)$ as a real vector space is just $\omega(MG_1, MG_2) 
\otimes \pi$ modulo the radical of $(,)_{\pi}$ (see ~\cite{li2},
	~\cite{theta}).
The main object of study in this paper is $\theta_s$.  \\
\\
If $\pi$ is in $\mc R_s(MG_1, MG_2)$ but not
in $\mc R(MG_1, \omega(MG_1, MG_2))$, our construction from ~\cite{theta} will
result in a vanishing $(,)_{\pi}$. Thus $\theta_s(MG_1, MG_2)(\pi)$ "vanishes".
In this case, $\theta_s=\theta$ trivially. The remaining question is
whether $(,)_{\pi} \neq 0 $ if $\pi \in \mc R(MG_1,
\omega(MG_1, MG_2))$. Conjecturally, $\theta_s(MG_1, MG_1)$ should agree with 
the restriction of $\theta(MG_1, MG_1)$ on $\mc R_s(MG_1, MG_2)$ (see 
~\cite{theta}, ~\cite{li1}). \\
\\
For $\pi$ a Hermitian representation, it can be easily shown that $(,)_{\pi}$ is 
an invariant Hermitian 
form on $\theta(MG_1, MG_2)(\pi)$ if $(,)_{\pi}$ does not vanish. This is a 
special case of Przebinda's result in (~\cite{pr2}). For $\pi$ unitary, we do 
not know whether $(,)_{\pi}$ must be positive semidefinite in general. 
Nevertheless, in ~\cite{unit}, we have proved the semi-positivity of $(,)_{\pi}$
under certain condition on the leading exponents of $\pi$ (see ~\cite{knapp}, ~\cite{wallach}). 
Fix a Cartan decomposition for $Sp_{2n}(\mathbb R)$ and $O(p,q)$. Fix the standard basis of $\f a$ for $Sp_{2n}(\mathbb R)$ and $O(p,q)$ (see 6.1). The leading exponents of an irreducible admissible representation are in the complex dual of the Lie algebra $\f a$ of $A$. 
\begin{thm}
Suppose $p+q \leq 2n+1$. Let $\pi$ be an irreducible unitary representation 
whose every leading exponent satisfies
\begin{equation}~\label{ss1}
\Re(v)-(\bold{n}-\frac{\bold{p+q}}{2})+\rho(O(p,q)) \preceq 0
\end{equation}
Then $(,)_{\pi}$ is positive semidefinite. Thus, $\theta_s(p,q;2n)(\pi)$ is 
either unitary or vanishes. 
\end{thm}
We denote the set of representations in $\Pi(MO(p,q))$ satisfying  (~\ref{ss1}) 
by $\mc R_{ss}(p,q;2n)$. The set 
$\mc R_s(MO(p,q),MSp_{2n}(\mb R))$ is written as $\mc R_s(p,q;2n)$ in short.
\begin{thm}
Suppose $n < p \leq q$. Let $\pi$ be an irreducible unitary representation whose 
every leading exponent satisfies
\begin{equation}~\label{ss2}
\Re(v)-(\frac{\bold{p+q}}{2}-\bold{n}-\bold{1})+\rho(Sp_{2n}(\mb R)) \preceq 0
\end{equation}
Then $(,)_{\pi}$ is positive semidefinite. Thus, either 
$\theta_s(p,q;2n)_s(\pi)$ is unitary or vanishes. 
\end{thm}
We denote the set of representations in $\Pi(MSp_{2n}(\mb R))$ satisfying  
(~\ref{ss2}) by $\mc R_{ss}(2n; p,q)$. The set 
$\mc R_s(MSp_{2n}(\mb R), MO(p,q))$ is written as $\mc R_s(2n; p,q)$ in short.\\
\subsection{Estimates on Leading Exponents and $L(p,n)$}
In this paper, we establish some estimates on the growth of the matrix 
coefficients of $\theta(p,q;2n)(\pi)$ and of $\theta(2n;p,q)(\pi)$ for $\pi$ in 
$\mc R_s(p,q;2n)$
and $\mc R_s(2n;p,q)$ respectively. We 
achieve this by studying the decaying of
the function
$$L(a, \phi)=\int_{b_1 \geq b_2 \ldots \geq b_p \geq 1} 
(\prod_{i=1,j=1}^{n,p}(a_i^2+ b_j^2)^{-\frac{1}{2}} )
\phi(b_1,b_2,\ldots,b_p) d b_1 d b_2 \ldots d b_p$$
as a function of $a \in \mb R^n$. In general, the decaying of $L(a, \phi)$ 
depends on the decaying of $\phi$. In section 5, we define a map 
$L(p,n)$ to describe this dependence. The map $L(p,n)$ is a continuous map from 
$$C(p)=\{\lambda \prec 0 \mid \lambda \in \mb R^p\}$$
 to 
$$C(n)=\{ \mu \prec 0 \mid \mu \in \mb R^n \}$$
Its algorithm is developed in Section $5$.
For some special vectors in $C(p)$, $L(p,n)$ is just a reordering plus an 
augmentation or truncation. In this paper, 
we prove
\begin{thm}.
Let $L(n,p)$ be defined as in Section $5$. Let $a(g_2)$ be the middle term of 
the $KA^+K$ decomposition of $g_2 \in Sp_{2n}(\mb R)$. Let $b(g_1)$ be the middle term 
of the $KA^+K$ decomposition of $g_1 \in O(p,q)$.
\begin{enumerate}
\item
Suppose that $\pi \in \mc R_s(p,q;2n)$. Suppose $\lambda \prec -2 \rho(O(p,q))+ 
\bold n $ and for every leading exponent $v$ of $\pi$, $\Re(v) \preceq \lambda$. 
Then
the matrix coefficients of $\theta_s(p,q;2n)(\pi)$ are weakly bounded by
$$a(g_2)^{L(p,n)(\lambda+2 \rho(O(p,q))-\bold n)-\bold{\frac{q-p}{2}}}.$$
\item Suppose that $\pi \in \mc R_s(2n;p,q)$. Suppose 
$\lambda \prec -2 \rho(Sp_{2n}(\mb R))+\frac{\bold{p+q}}{2}$ and for every 
leading exponent $v$ of $\pi$, $\Re(v) \preceq \lambda$.
 Then
the matrix coefficients of $\theta_s(2n;p,q)(\pi)$ are weakly bounded by
$$b(g_1)^{L(n,p)(\lambda+2 \rho(Sp_{2n}(\mb R))-\frac{\bold{p+q}}{2})}$$
\end{enumerate}
\end{thm}
The definition of weakly boundedness is given in Section $3$.
\subsection{Quantum Induction}
The idea of composing two theta correspondences to obtain "new" representations 
has been known for years. For example, one can compose $\theta(p,q;2n)$ with 
$\theta(2n;p^{\p},q^{\p})$. The nature of 
$\theta(2n;p^{\p},q^{\p}) \theta(p,q;2n)(\pi)$ seems to be inaccessible except 
for the cases of stable ranges. In this paper, we treat a somewhat more 
accessible object, namely, 
$$\theta_s(2n;p^{\p},q^{\p}) \theta_s(p,q;2n)(\pi).$$
Our construction is done through the studies of the Hermitian form $(,)_{\pi}$. 
Due to the unitarity theorems we proved in (~\cite{unit}), under restrictions as 
specificed in Equations (~\ref{ss1}) and (~\ref{ss2}),
quantum induction preserves unitarity. 
Our main result can be stated as follows
\begin{thm}[Main Theorem]
\begin{itemize}
\item Suppose
\begin{enumerate} 
\item $q^{\p} \geq p^{\p} > n$;
\item $p^{\p}+q^{\p}-2n \geq 2n-(p+q)+2 \geq 1$;
\item $p+q =p^{\p}+q^{\p} \qquad (\mod \ 2)$.
\end{enumerate}
Let $\pi$ be an irreducible unitary representation in $\mc R_{ss}(p,q;2n)$. 
Suppose that $(,)_{\pi}$ does not vanish. Then
\begin{enumerate}
\item  $\theta_s(p,q;2n)(\pi)$ is unitary.
\item  $\theta_s(p,q;2n)(\pi) \in \mc R_{ss}(2n;p^{\p},q^{\p})$.
\item 
$\theta_s(2n; p^{\p},q^{\p})\theta_s(p,q;2n)(\pi)$ is either an irreducible
 unitary representation or the NULL representation.
\end{enumerate}
\item Suppose 
\begin{enumerate}
\item $2n^{\p}-p-q+2 \geq p+q-2n$;
\item  $n < p \leq q$.
\end{enumerate}
Let $\pi$ be a unitary representation in $\mc R_{ss}(p,q;2n)$. Suppose 
$(,)_{\pi}$ does not vanish. Then 
\begin{enumerate}
\item $\theta_s(2n;p,q)(\pi)$ is unitary.
\item $\theta_s(2n;p,q)(\pi) \in \mc R_{ss}(p,q;2n^{\p})$.
\item $\theta_s(p,q;2n^{\p})\theta_s(2n;p,q)(\pi)$ is either an irreducible 
unitary
representation or the NULL representation.
\end{enumerate}
\end{itemize}
\end{thm}
The purpose of assuming $\pi \in \mc R_{ss}$ is to guarantee the unitarity of $
Q(*)(\pi)$. In fact, for any $\pi$, the condition on the sizes of related dual 
pairs can be computed easily to define nonunitary quantum induction.
In general, the underlying Hilbert space of the induced representation is 
"invisible" under quantum induction except for certain limit cases where quantum 
induction becomes unitary parabolic induction (see Section 6 and ~\cite{quan}). 
\begin{conj} Suppose $\pi$ is a unitary representation in $\mc R_{ss}$.
\begin{itemize}
\item The quantum induction $Q(p,q;2n;p^{\p},q^{\p})(\pi)$ for $2n-p-
q+2=p^{\p}+q^{\p}-2n$ can be obtained via unitarity-preserving parabolic 
induction and cohomological induction from $\pi$.
\item The quantum induction $Q(2n;p,q;2n^{\p})(\pi)$ for $p+q-2n-2=2n^{\p}-p-q$ 
can be obtained as a subfactor via unitarity-preserving parabolic induction  
from $\pi$.
\end{itemize}
\end{conj}
For the cases $p+q=2n+1=p^{\p}+q^{\p}$ and $p+q=2n+1=2n^{\p}+1$, by a Theorem of 
Adams-Barbasch, $Q$ is either the identity map or vanishes (~\cite{ab}). Our 
conjecture holds trivially, i.e., no induction is needed. For the case 
$p+q+p^{\p}+q^{\p}=4n+2$ and $p-p^{\p}=q-q^{\p}$, our result in Section 6 gives 
some indication that $Q(p,q;2n;p^{\p},q^{\p})(\pi)$ can be obtained from
$$Ind_{SO_0(p,q) GL_0(p^{\p}-p) N}^{SO_0(p^{\p},q^{\p})}(\pi \otimes 1).$$

Let me make one remark regarding the nonvanishing of $(,)_{\pi}$.  In 
~\cite{non1} we prove
\begin{thm}[~\cite{non1}]
Suppose $p+q \leq 2n+1$. Let $\pi \in \mc R_s(p,q;2n)$. Then
at least one of
$$(,)_{\pi}, (,)_{\pi \otimes \det}$$
does not vanish.
\end{thm}
For $\pi \in \mc R_s(2n;p,q)$, the nonvanishing of $(,)_{\pi}$ is hard to detect 
since it depends on $p,q$ (~\cite{ab}, ~\cite{thesis}, ~\cite{moeglin}).
A result of Jian-Shu Li says that  $(,)_{\pi}$ does not vanish if $p,q \geq 2n$.
We are not aware of any more general nonvanishing theorems. \\
\\
Finally, concerning the associated varieties, Przebinda shows that the 
associated varieties behaves reasonably well under theta correspondence under 
certain strong hypothesis (~\cite{pr1}). We conjecture that quantum induction 
induces an induction on associated varieties and wave front sets. The exact 
description of the associated variety under quantum induction can be predicted 
based on ~\cite{pdk}.
\begin{conj} 
\begin{itemize}
\item Under the same assumptions from the main theorem, let $\pi$ be a unitary 
representation in $\mc R_{ss}(p,q;2n)$. Let $\mc O_{\bold d}$ be the associated 
variety of $\pi$ with
$\bold d$ a partition (see Ch 5, ~\cite{cm}). Let
$\mc O_{\bold f}$ be the associated variety of $Q(p,q;2n;p^{\p},q^{\p})(\pi) 
\neq 0$.
Then $\bold f^t=(p^{\p}+q^{\p}-2n, 2n-p-q, \bold d^t)$.
\item Under the same assumptions from the main theorem, let $\pi$ be a unitary 
representation in $\mc R_{ss}(2n;p,q)$. Let $\mc O_{\bold d}$ be the associated 
variety of $\pi$ with $\bold d$ a partition. Let $\mc O_{\bold f}$ be the 
associated variety of $Q(2n;p,q;2n^{\p})(\pi) \neq 0$. Then
$\bold f^t=(2n^{\p}-p-q, p+q-2n, \bold d^t)$.
\end{itemize}
\end{conj}
We remark that our situation is different from the situation treated in 
~\cite{pr1} with some overlaps.  The description of the wave front set under 
quantum induction can be predicted based on ~\cite{pan}.

\section{Theta Correspondence}
Let $(O(p,q), Sp_{2n}(\mb R))$ be a reductive dual pair in $Sp_{2n(p+q)}(\mb 
R)$. Let 
$$j: Mp_{2n(p+q)}(\mb R) \rightarrow Sp_{2n(p+q)}(\mb R)$$ be the double 
covering. Let $\{1, \epsilon\}=j^{-1}(1)$.
Let $MO(p,q)=j^{-1}(O(p,q))$ and $MSp_{2n}(\mb R)=j^{-1}(Sp_{2n}(\mb R))$. Fix a 
maximal compact subgroup $U$ of $Sp_{2n(p+q)}(\mb R)$ such that
$$U \cap Sp_{2n}(\mb R) \cong U(n), \qquad U \cap O(p,q) \cong O(p) \times 
O(q).$$
Then $MU$ is a maximal compact subgroup of $Mp_{2n(p+q)}(\mb R)$. Let 
$\omega(p,q;2n)$ be the oscillator representation of $Mp_{2n(p+q)}(\mb R)$.
The representation $\omega(p,q;2n)$ or sometimes $\omega(2n;p,q)$ is regarded as 
an admissible representation of $Mp_{2n(p+q)}(\mb R)$ equipped with a fixed dual 
pair $(O(p,q), Sp_{2n}(\mb R))$.
Let $\mc P$ be the Harish-Chandra module.
Then $\omega(p,q;n)$ can be restricted to $MO(p,q)$ and $MSp_{2n}(\mb R)$.
Howe's theorem states that there is a one to one correspondence 
$$\theta(p,q;2n): \mc R(MO(p,q), \omega(p,q; 2n)) \rightarrow \mc R(MSp_{2n}(\mb 
R), \omega(p,q;2n)).$$
\subsection{$MO(p,q)$ and $MSp_{2n}(\mb R)$}
The groups $MO(p,q)$ and $MSp_{2n}(\mb R)$ are double covers of $O(p,q)$ and
 $Sp_{2n}(\mb R)$. Depending on the parameter $n$, $p$ and $q$, they may be 
quite different.
\begin{lem}~\label{msp2n}
\begin{enumerate}
\item
If $p+q$ is odd, then the double cover $MSp_{2n}(\mb R)$ does not split. It is
the metaplectic group $Mp_{2n}(\mb R)$. The representations in 
$\mc R(Mp_{2n}(\mb R), \omega(p,q; 2n))$ are
genuine representation of $Mp_{2n}(\mb R)$. 
\item If $p+q$ is even, then the double cover
$MSp_{2n}(\mb R)$ splits. It is the product of $Sp_{2n}(\mb R)$ and $\{1, 
\epsilon \}$.
The representations in $\mc R(MSp_{2n}(\mb R), \omega(p,q;2n))$ can be 
identified with
representations of $Sp_{2n}(\mb R)$ by tensoring the 
nontrivial character of $\{1, \epsilon\}$. 
\item In both cases, any representation in
$$\mc R(MSp_{2n}(\mb R), \omega(p,q; 2n))$$
can be identified with a representation of $Mp_{2n}(\mb R)$. In the former case, 
a
genuine representation, and in the latter case, a nongenuine representation.
\end{enumerate}
\end{lem}
We do not know the earliest reference. The details can be worked out easily and 
can be found in ~\cite{ab}.
\begin{lem}~\label{mopq}
\begin{enumerate}
\item As a group,
$$MO(p,q) \cong \{ (\xi, g) \mid g \in O(p,q), \xi^2=\det g^n \}$$
\item
$\xi$ is a character of $MO(p,q)$. Any representations in $\mc R(MO(p,q), 
\omega(p,q;2n))$
can be identified with representations of $O(p,q)$ by tensoring $\xi$.
\item $MSO(p,q)$ can be identified as group product
$$SO(p,q) \times \{1, \epsilon\}.$$
\item If $n$ is even, $MO(p,q) \cong O(p,q) \times \{1, \epsilon\}$.
\end{enumerate}
\end{lem}
The details can be found in ~\cite{ab} or ~\cite{unit}.
We must keep in mind that for $p+q$ odd,
$$\mc R(MSp_{2n}(\mb R), \omega(p,q;2n)) \subset \Pi_{genuine}(Mp_{2n}(\mb R))$$
 and for $p+q$ even
$$\mc R(MSp_{2n}(\mb R), \omega(p,q;2n)) \subset \Pi(Sp_{2n}(\mb R)).$$

\subsection{Averaging Integral $(,)_{\pi}$}
Let $O(p,q)$ be the orthogonal group preserving the symmetric form defined by
$$I_{p,q}=\left( \begin{array}{clcr} 0_p & 0 & I_p \\ 0 & I_{q-p} & 0 \\ I_p & 0 
& 0_p
 \end{array} \right).$$ 
Fix a Cartan decomposition with
$$A=\{ \diag(a_1,a_2, \ldots, a_p,\overbrace{1,\ldots,1}^{q-p},a_1^{-1}, a_2^{-1}, \ldots, a_p^{-1}) \mid a_i >0 \}$$
and a positive Weyl chamber
$$A^+=\{ \diag(a_1, a_2 \ldots, a_p,\overbrace{1,\ldots,1}^{q-p},a_1^{-1}, a_2^{-1}, \ldots, a_p^{-1}) \mid a_1 \geq a_2 \geq \ldots \geq a_p \geq 1  \}.$$
The half sum of the positive restricted roots of $O(p,q)$
$$\rho(O(p,q))=\overbrace{(\frac{p+q-2}{2}, \frac{p+q-4}{2}, \ldots, \frac{q-
p}{2})}^{p}.$$

Let $Sp_{2n}(\mb R)$ be the symplectic group that preserves the skew-symmetric 
form defined
by 
$$W_n=\arr{0_n & -I_n \\ I_n & 0_n }$$
Let $K$ be the intersection of $\xin$ with the orthogonal group $O(2n)$ which 
preserves the Euclidean
inner product on $\mb R^{2n}$. Let
$$A=\{a= \diag(a_1, a_2, \ldots, a_n, a_1^{-1}, \ldots, a_n^{-1}) \mid a_i > 
0 \}$$
$$A^+=\{a= \diag(a_1, a_2, \ldots, a_n, a_1^{-1}, \ldots, a_n^{-1}) \mid a_1 
\geq a_2 \geq \ldots
\geq a_n \geq 1\}.$$
The half sum of the positive restricted roots of $Sp_{2n}(\mb R)$
$$\rho(Sp_{2n}(\mb R))=\overbrace{(n,n-1, \ldots, 1)}^n.$$

For each irreducible admissible representation of a semisimple group $G$ of real 
rank $r$, 
there are number of $r$-dimensional complex vectors in $\f a^*$ called leading exponents 
attached to it.
Leading exponents are the main data used to produce the Langlands classification 
(see ~\cite{langlands}
and ~\cite{knapp}).
\begin{defn}
An irreducible representation $\pi$ of $O(p,q)$ is said to be
in the {\it semistable range} of
$\theta(p,q;2n)$
if and only if each leading exponent $v$ of $\pi$ satisfies
\begin{equation}
\sum_{i=1}^{j} \Re(v_i)+(p+q-2i)-n < 0 \qquad (\forall \ \ j \in [1,p])
\end{equation}
i.e.,
$$\Re(v)-\bold n+ 2 \rho(O(p,q)) \prec 0.$$
An irreducible representation $\pi$ of $Mp_{2n}(\mb R)$ is said to be
in the semistable range of $\theta(2n;p,q)$ if and only if every leading 
exponent
$v$ of $\pi$ satisfies
\begin{equation}
\sum_{i=1}^k \Re (v_i)-\frac{p+q}{2}+2n+2-2j < 0 \qquad (\forall \ \ k \in 
[1,n])
\end{equation}
i.e.,
$$\Re(v) -\frac{\bold{p+q}}{2}+ 2 \rho(Sp_{2n}(\mb R)) \prec 0.$$
\end{defn}
If $W$ is a complex linear space, we use a superscript $W^c$ to denote $W$ 
equipped with
the conjugate 
complex linear structure.
Let $\pi \in \mc R_s(MG_1,MG_2)$. We define a complex linear
pairing
$$ (\mc P^c \otimes \pi, \mc P \otimes \pi^c ) \rightarrow \mb C$$
as follows: for $\phi \in \mc P, \psi  \in \mc P^c, v \in \pi^c, u \in \pi$,
$$( \phi \otimes v , \psi \otimes u)_{\pi}
=\int_{MO(p,q)} (\phi, \omega(g) \psi) ( \pi(g)u, v)
d g$$
If $\pi$ is unitary, $(,)_{\pi}$ is an invariant Hermitian form with respect to 
the action of $MG_2$.
 
\begin{thm}(see ~\cite{theta})
Suppose $(\pi, V)$ is a unitary representation in the semistable range of 
$\theta(MG_1,MG_2)$.
 Then $(,)_{\pi}$ is well-defined.
Suppose $\mc R_{\pi}$ is the radical of $(,)_{\pi}$ with
respect to $\mc P \otimes V^c$. If $(,)_{\pi}$  does not vanish,
then
\begin{itemize}
\item $\pi$ occurs in $\mc R(MG_1, \omega(MG_1,MG_2))$;
\item $\mc P \otimes V^c / \mc R_{\pi}$ is irreducible;
\item $\mc P \otimes V^c / \mc R_{\pi}$  is isomorphic to 
$\theta(MG_1, MG_2)(\pi)$.
\item $\theta_s(MG_1,MG_2)(\pi)$ is a Hermitian representation of $MG_2$.
\end{itemize}
\end{thm}
Thus the Harish-Chandra module of $\theta_s(MG_1, MG_2)(\pi)$ can be defined as
$\mc P \otimes V^c/\mc R_{\pi}$.
\subsection{Oscillator Representation}
The oscillator representation, also known as the Segal-Shale-Weil representation, 
is a unitary representation of the metaplectic group $Mp$. The construction of the oscillator representation can be found in
the papers of Segal, Shale and Weil (~\cite{shale},~\cite{segal}, ~\cite{weil}). In this section, we give a basic estimate of the matrix coefficients of the oscillator representation. Proof of Theorem 3.3.1 can also be found in ~\cite{howe82} (Prop. 8.1).\\
\\
Let $g \in \xin$.
Let $a(g)$ be the midterm of the $KAK$ decomposition of $g$ such that $a \in 
A^+$. Let $H(g)=\log a(g)$.
Then
$$H(g)=\diag(H_1(g), H_2(g), \ldots , H_n(g), -H_1(g), \ldots, -H_n(g))$$
is in the Weyl chamber $\f a^+$.\\
\\
Let $Mp_{2n}(\mb R)$ be the double covering of $\xin$. The midterm of the $KAK$ 
decomposition
of $Mp_{2n}(\mb R)$ remains the same. Let $(\omega_n, L^2(\mb R^n))$ be the 
Schr\"odinger
model of the oscillator representation of $Mp_{2n}(\mb R)$ as in ~\cite{theta}. 
Let
$$\mu(x)=\exp (-\frac{1}{2}(x_1^2+x_2^2 + \ldots + x_n^2))$$
be the Gaussian function.
The Harish-Chandra module $\mc P_n$ are the polynomial functions multiplied by 
the Gaussian function
as verified in ~\cite{theta}. We write
$$x^{\alpha}=\prod_1^n x_i^{\alpha_i}.$$
Harish-Chandra's theory says that the $Mp_{2n}(\mb R)$ action on $\mc P_n$ can 
be controlled by
the $A$ action on fixed $K$-types of $\omega_n$.
\begin{thm}~\label{metaplectic} For any $a \in A$, we have
$$(\omega_n(a) x^{\alpha} \mu(x), x^{\beta} \mu(x))=c_{\alpha,\beta} 
\prod_{i=1}^n 
a_i^{\alpha_i+\frac{1}{2}}(1+a_i^2)^{-\frac{\alpha_i+\beta_i+1}{2}}$$
In addition,
$$|(\omega_n(a) x^{\alpha} \mu(x), x^{\beta} \mu(x))| \leq c \prod_{i=1}^n 
(a_i+a_i^{-1})^{-\frac{1}{2}}$$
In general, for every $\phi, \psi \in \mc P_n$, we have
$$|(\omega_n(g) \phi, \psi) | \leq c \prod_{i=1}^n (a_i(g)+a_i^{-1}(g))^{-
\frac{1}{2}}$$
\end{thm}
The proof for the first statement can be found in ~\cite{theta}. We observe that
\begin{equation}
\begin{split}
  & |(\omega_n(a) x^{\alpha} \mu(x), x^{\beta} \mu(x)) |\\
= & |c_{\alpha,\beta} \prod_{i=1}^n 
a_i^{\alpha_i+\frac{1}{2}}(1+a_i^2)^{-\frac{\alpha_i+\beta_i+1}{2}} | \\
= & |c_{\alpha,\beta} \prod_{i=1}^n 
(a_i+a_i^{-1})^{-\frac{1}{2}}(1+a_i^2)^{-\frac{\beta_i}{2}}
(1+a_i^{-2})^{-\frac{\alpha_i}{2}} | \\
\leq & c_{\alpha,\beta} \prod_{i-1}^n (a_i +a_i^{-1})^{-\frac{1}{2}}.
\end{split}
\end{equation} 
The second statement is proved. The third statement follows immediately from
$K$-finiteness of $\phi$ and $\psi$. \\
\\
The estimations on the right hand side are invariant under Weyl group action, 
thus do not depend on the choices
of the Weyl chamber $\f a^+$. 
\subsection{Growth of Matrix Coefficients}
\begin{defn} Suppose $X$ is a Borel measure space equipped with a norm $\|.\|$ 
such that 
\begin{itemize}
\item $\|x \| \geq 0$ for all $x \in X$;
\item the set $\{ \|x \| \leq r \}$ is compact.
\end{itemize}
Let $f(x)$ and $\phi(x)$ be continuous functions defined over $X$.
Suppose $\phi(x)$ approaches $0$ as $\|x \| \rightarrow \infty$.
A function $f(x)$ is said to be weakly bounded by the function $\phi(x)$ if 
there exists a $\delta_0>0$ such that
for every $\delta_0 > \delta >0$, there exists a $C >0$ depending on $\delta$ 
such that
$$|f(x)| \leq C  \phi(x)^{1-\delta} \qquad (\forall \ x \in X)$$
\end{defn}
The typical case is when $f(x)$ does not decay as fast as $\phi(x)$ but faster 
than $\phi(x)^{1-\delta}$.\\
\\
Let $\pi$ be an irreducible representation of a reductive group $G$. Let $K$ 
be a maximal compact subgroup of $G$. We adopt the notation from Chapter VIII in 
~\cite{knapp}. We equip $G$ with a norm 
$$ g \rightarrow \|\log(a(g))\|=(\log a(g), \log a(g))^{\frac{1}{2}}$$
where $(,)$ is a real $\f g$-invariant symmetric form whose restriction on $\f a$ is 
positive definite. \\
\\
{\bf Example}: An irreducible representation $\pi$ of a reductive group $G$ is 
tempered 
if and only if
its matrix coefficients are weakly bounded by 
$$a(g)^{ -\rho}$$
where $\rho$ is the half sum of positive restricted roots and $a(g)$ is the mid 
term of
the $KAK$ decomposition with $a(g)$ in the positive Weyl chamber $A^+$ (see 
~\cite{knapp}). 
\begin{thm}~\label{equ}
Let $\pi$ be an irreducible unitary representation of $G$. Let $\lambda \prec 
0$.
The following are equivalent.
\begin{enumerate}
\item Every leading exponent $v$ of $\pi$ has $\Re(v) \preceq \lambda$.
\item There is an integer $q \geq 0$ such that every $K$-finite matrix 
coefficient is bounded by a multiple of $(1+\|\log a(g)\|)^q \exp(\lambda(\log 
a(g)))$.
\item Every $K$-finite matrix coefficient $\phi(g)$ of $\pi$ is bounded by $C 
a(g)^{\lambda+ \delta}$ for any $\delta \succ 0$.
\item Every $K$-finite matrix coefficient of $\pi$ is weakly bounded by 
$a(g)^{\lambda}$.
\end{enumerate}
\end{thm}
See Chapter VIII.8,13 ~\cite{knapp} or Chapter 4.3 ~\cite{wallach} for details. 
The first three statements are equivalent without assuming the unitarity of 
$\pi$ and $\lambda \prec 0$.

\section{Twisted Integral}
Let $A^+=\{a_1 \geq a_2 \ldots \geq 1 \}$.
In this section, we will study the following integrals
$$L(a,\lambda)=\int_{B^+} \prod_{i=1}^p (\prod_{k=1}^n (a_k^2+b_i^2)^{-
\frac{1}{2}})
 b_i^{\lambda_i} d b_i$$
and
$$L(a, \phi)=\int_{b_1 \geq b_2 \ldots \geq b_p \geq 1} \prod_{i,j}(a_i^2+ 
b_j^2)^{-\frac{1}{2}} 
\phi(b_1,b_2,\ldots,b_p) d b_1 d b_2 \ldots d b_p .$$
The domain of $a$ will always be $A^+$ unless stated otherwise.
We are interested in the growth of $L(a,\phi)$ as $a$ goes to infinity. 
Variables and parameters are assumed to be real in this section.
\subsection{Single Variable Case $a \geq 1$}
\begin{lem} Suppose that $a \geq 1$.
The integral
$$L(a, \lambda)=\int_{b \geq 1} (a^2+b^2)^{-\frac{1}{2}} b^{\lambda} {d b}$$
converges if and only if $\lambda < 0$. In addition, $L(a, \lambda)$ is weakly
 bounded by $a^{\lambda}$ if $-1 \leq \lambda<0$ and is  bounded by a multiple 
of
$a^{-1}$ if $\lambda < -1$.
\end{lem}
Proof: From classical analysis, the integral
$$\int_{b \geq 1}  b^{-1+\lambda} {d b}$$
 converges if and only if
$\lambda<0$. 
For a fixed $a$ and any $b>1$, $b^2 \leq a^2+b^2 \leq (1+a^2) b^2$. Hence
$$\int_{b \geq 1} b^{-1} b^{\lambda} d b \geq \int_{b \geq 1} (a^2+b^2)^{-
\frac{1}{2}} b^{\lambda} {d b} \geq \int_{b \geq 1} (1+a^2)^{-\frac{1}{2}} b^{-
1} b^{\lambda} d b$$
Hence,
$ L(a, \lambda)$ converges if and only if $\lambda<0$. \\
\\
For $a \geq 1$,
\begin{equation}
\begin{split}
L(a, \lambda) = & \int_{b \geq 1} (a^2+b^2)^{-\frac{1}{2}} b^{\lambda} d b \\
= & \int_{ab \geq 1} (a^2+a^2 b^2)^{-\frac{1}{2}} a^{\lambda+1} b^{\lambda} d b \\
= & a^{\lambda} \int_{b \geq a^{-1}} (1+b^2)^{-\frac{1}{2}} b^{\lambda} d b \\
= & a^{\lambda} \int_{b \geq 1} (1+b^2)^{-\frac{1}{2}} b^{\lambda} d b+ a^{\lambda} \int_{a^{-1}}^{1} (1+b^2)^{-\frac{1}{2}} b^{\lambda} d b
\end{split}
\end{equation}
For $a \geq 1$ and $a^{-1} \leq b \leq 1$ and $\lambda \neq -1$,
$$\frac{1}{\sqrt{2}} b^{\lambda} \leq (1+b^2)^{-\frac{1}{2}} b^{\lambda} \leq b^{\lambda}.$$
 Taking $\int_{a^{-1}}^1 d b $, we obtain
$$\frac{1}{\sqrt{2}(\lambda+1)}(a^{\lambda}-a^{-1}) \leq a^{\lambda} \int_{a^{-1}}^{1} (1+b^2)^{-\frac{1}{2}} b ^{\lambda} d b
\leq \frac{1}{\lambda+1}(a^{\lambda}-a^{-1}) .$$
Therefore, for $-1 < \lambda < 0$, $L(a, \lambda)$ is bounded by a multiple of $a^{\lambda}$; for $\lambda < -1$, $L(a, \lambda)$ is bounded by a multiple of $a^{-1}$. For $\lambda=-1$,
$$\frac{1}{\sqrt{2}}a^{-1} \ln a \leq a^{-1} \int_{a^{-1}}^{1} (1+b^2)^{-\frac{1}{2}} b ^{-1} d b
\leq a^{-1} \ln a.$$
Therefore, $L(a, -1)$ is weakly bounded by $a^{-1}$. Q.E.D. \\
\\
\begin{lem}~\label{delta}
Suppose $\lambda_0 <0$. Suppose $f(a)$ is weakly bounded by $a^{\lambda}$ for 
any $0>\lambda > \lambda_0$. Then $f(a)$ is weakly bounded by $a^{\lambda_0}$.
\end{lem}
Combining these two lemmas, we obtain
\begin{thm} Suppose that $a \geq 1$. Suppose $\phi(b)$ is weakly bounded by 
$b^{\lambda}$ for some $\lambda <0$.
Then the integral
$$L(a,  \phi(b))=\int_{b \geq 1} (a^2+b^2)^{-\frac{1}{2}} \phi(b) d b$$
converges. In addition, $L(a, \phi)$ is weakly
 bounded by $a^{\lambda}$ if $-1 \leq \lambda$ and is bounded by a multiple of
$a^{-1}$ if $\lambda < -1$.
\end{thm}
In conclusion, the growth rate of $L(a, \phi(b))$ is a
"truncation" of the growth rate of $\phi(b)$.

\subsection{Multivariate $b$}
Let $\lambda=(\lambda_1,\lambda_2,\ldots,\lambda_p)$.
Let $B^+=\{b_1 \geq b_2 \geq \ldots \geq b_p \geq 1\}$.
Let us consider
$$L(a, \lambda)=\int_{B^+} \prod_{i=1}^p (a^2+b_i^2)^{-\frac{1}{2}} 
b_i^{\lambda_i} d b_i $$
First, we observe that
$$a^2+b_i^2 \geq a^{2 \eta_i} b_i^{2-2 \eta_i}$$
for any $\eta_i \in [0,1]$. The $\eta_i$ is to be determined later.
We obtain
\begin{equation}
\begin{split}
L(a, \lambda) \leq & \int_{B^+} \prod_{i=1}^p a^{-\eta_i} b_i^{-
1+\eta_i+\lambda_i} d b_i \\
= & a^{\sum_{i=1}^p -\eta_i} \int_{B^+} \prod_{i=1}^p b_i^{-1+\eta_i+\lambda_i} 
d b_i \\
\end{split}
\end{equation}
Secondly, we change the coordinates and let
$$r_i=\frac{b_i}{b_{i+1}} \qquad (i=1,\ldots, p-1)$$
$$r_p=b_p.$$
Then
$$b_i=\prod_{j=i}^{p} r_j \qquad (i=1,\ldots, p).$$
In addition, $B^+$ is transformed into $[1, \infty)^p$.
The differential
$$ \prod_{i=1}^p d b_i=\prod_{i=1}^p (\prod_{j=i}^p r_j) \frac{d 
r_i}{r_i}=\prod_{i=1}^p b_i
\frac{ dr_i}{r_i}.$$
We obtain
\begin{equation}
\begin{split}
L(a, \lambda) & \leq a^{-\sum_{i=1}^p \eta_i}
\int_{[1, \infty)^p} \prod_{i=1}^p b_i^{ \eta_i+\lambda_i} \frac{ dr_i}{r_i} \\
 & =a^{-\sum_{i=1}^p \eta_i}
\int_{[1, \infty)^p} \prod_{i=1}^p (\prod_{j=i}^p r_j^{ \eta_i+\lambda_i}) 
\frac{ d r_i}{r_i} \\
&=a^{-\sum_{i=1}^p \eta_i}
\int_{[1, \infty)^p} \prod_{j=1}^p r_j^{\sum_{i=1}^j \eta_i+\lambda_i} \frac{d 
r_j}{r_j}.
\end{split}
\end{equation}
This integral converges if
$$\sum_{i=1}^j \eta_i+\lambda_i < 0 \qquad (\forall \, \, j).$$
\begin{thm} Suppose $ a \geq 1$. If $\lambda \prec 0$, then
$L(a, \lambda)$ converges. Furthermore, $L(a, \lambda)$ is bounded by a multiple 
of
$$a^{\sum_{i=1}^p \eta_i}$$
 with any $\eta_i$ satisfying the condition
$$\{ 0 \leq \eta_j \leq 1 , \sum_{i=1}^j 
\eta_i + \sum_{i=1}^j 
\lambda_i <0  \qquad (j=1, \ldots p) \}.$$
\end{thm}
The condition
$$\sum_{i=1}^j 
\eta_i + \sum_{i=1}^j 
\lambda_i <0  \qquad (j=1, \ldots p) $$
can be restated as $\eta+\lambda \prec 0$.
Combined with Lemma ~\ref{delta}, we have
\begin{thm} Suppose $\phi(b_1, b_2, \ldots, b_p)$ on $B^+$ is weakly bounded by 
$b^{\lambda}$
for some $\lambda \prec 0$. Then the function
$$L(a, \phi)= \int_{B^+} (\prod_{i=1}^p (a^2+b_i^2)^{-\frac{1}{2}}) \phi(b) d b_1 
\ldots d b_p$$
is weakly bounded by $a^{-\mu}$ with
$$\mu=\max \{\sum_{i=1}^p \eta_i \mid 0 \leq \eta_j \leq 1 , \lambda+\eta 
\preceq  0 \}.$$
\end{thm}
We point out the second ingredient needed to carry out estimations on 
$L(a,\phi)$, namely, the coordinate transform from $b$ to $r$.

\subsection{Multivariate $a \in [1,\infty)^n$}
This case is more complicated since the function
$L(a, \phi)$ is no longer of single variable. Our 
result here is weaker than
the results for single variable $a$. \\
\\
First we consider
$$L(a,\lambda)=\int_{B^+} \prod_{i=1}^p (\prod_{k=1}^n (a_k^2+b_i^2)^{-
\frac{1}{2}})
 b_i^{\lambda_i} d b_i $$
We again set the parameters $\eta_{k,i}$ to be in $[0,1]$. We have
$$a_k^2+b_i^2 \geq a_k^{2 \eta_{k,i}} b_i^{2-2\eta_{k,i}}$$
Therefore, we obtain
\begin{equation}
\begin{split}
 & L(a, \lambda) \\
\leq & \int_{B^+} \prod_{i=1}^p (\prod_{k=1}^n a_k^{-\eta_{k,i}}
b_i^{-1+\eta_{k,i}}) b_i^{\lambda_i} d b_i \\
= & \prod_{k=1}^n a_k^{-\sum_{i=1}^p \eta_{k,i}}
\int_{B^+} \prod_{i=1}^p b_i^{\lambda_i-n+\sum_{k=1}^n \eta_{k,i}} d b_i
\end{split}
\end{equation}
Now we change the coordinates $b$ into $r$. We obtain
\begin{equation}
\begin{split}
& L(a, \lambda) \\
\leq & \prod_{k=1}^n a_k^{-\sum_{i=1}^p \eta_{k,i}}
\int_{[1, \infty)^p} \prod_{i=1}^p 
(\prod_{j=i}^p r_j^{\lambda_i-n+\sum_{k=1}^n \eta_{k,i}} 
\prod_{j=i}^p r_j) \frac{ d r_i}{r_i} \\
= & \prod_{k=1}^n a_k^{-\sum_{i=1}^p \eta_{k,i}}
\int_{[1, \infty)^p} \prod_{i=1}^p 
(\prod_{j=i}^p r_j^{\lambda_i-n+1+\sum_{k=1}^n \eta_{k,i}} )
\frac{ d r_i}{r_i} \\
=& \prod_{k=1}^n a_k^{-\sum_{i=1}^p \eta_{k,i}}
\int_{[1, \infty)^p} \prod_{j=1}^p r_j^{\sum_{i=1}^j (\lambda_i-n+1+\sum_{k=1}^n 
\eta_{k,i})}
\frac{ d r_j}{ r_j}
\end{split}
\end{equation}
This integral converges if 
$$\sum_{i=1}^j (\lambda_i-n+1+\sum_{k=1}^n \eta_{k,i}) <0 \qquad (\forall \ \ 1 
\leq j \leq p).$$
Since $\eta_{k,i} \in [0,1]$,  we obtain the following theorem.
\begin{thm} Suppose $a \in [1,\infty)^p$.
The integral $L(a, \lambda)$ converges if 
$$\sum_{i=1}^j \lambda_i-n+1 <0$$
for every integer $1 \leq j \leq p$. In this situation 
$L(a,\lambda)$ is bounded by a multiple of
$$a^{-\mu}=\prod_{k=1}^n a_k^{-\mu_k}$$
where $\mu_k= \sum_{i=1}^p \eta_{k,i}$ and $\{ \eta_{k,i} \}$ satisfy
$$\eta_{k,i} \in [0,1] \ \ \forall k,i$$
\begin{equation}
\sum_{i=1}^j (\lambda_i-n+1+\sum_{k=1}^n \eta_{k,i}) < 0 \qquad \forall \ \ j.
\end{equation}
\end{thm}
Similarly, we obtain
\begin{thm}~\label{tw1}
Suppose $a \in [1,\infty)^p$. Suppose $\phi(b)b^{-\bold{n}+{1}}$ on $B^+$
is bounded  by $b^{\lambda}$ with $\lambda \prec 0$. 
Then the integral $L(a, \phi)$ converges. Furthermore,
$L(a,\phi)$ is bounded by a multiple of
$$a^{-\mu}=\prod_{k=1}^n a_k^{-\mu_k}$$
where $\mu_k= \sum_{i=1}^p \eta_{k,i}$ and $\{ \eta_{k,i} \}$ satisfy
$$\eta_{k,i} \in [0,1] \ \ \forall k,i$$
\begin{equation}~\label{tw11}
\sum_{i=1}^j (\lambda_i+\sum_{k=1}^n \eta_{k,i}) < 0 \qquad  \forall \ \  j.
\end{equation}
\end{thm}
\section{Algorithm and Examples}
Suppose $\lambda \prec 0$.  We are interested in finding the "maximal" $\eta$ where
$$\mu_k= \sum_{i=1}^p \eta_{k,i}$$
with $\eta_{k,i}$ satisfying
$$\eta_{k,i} \in [0,1] \ \ \forall k,i$$
\begin{equation}~\label{basic}
\sum_{i=1}^j (\lambda_i+\sum_{k=1}^n \eta_{k,i}) < 0 \qquad \forall \ \ j.
\end{equation}
\subsection{A Theorem for $a \in [1, \infty)^n$}
Write (~\ref{basic}) as
\begin{equation}~\label{basic1}
\sum_{i=1}^j (\sum_{k=1}^n \eta_{k,i}) < -\sum_{i=1}^j \lambda_i \qquad \forall 
\ \ j.
\end{equation}
First of all, since $\eta_{k,i} \geq 0$, the sequence
$$\{\sum_{i=1}^j \sum_{k=1}^n \eta_{k,i} \mid j \in [1, p] \}$$
is increasing. However, the sequence
$$\{-\sum_{i=1}^j \lambda_i \mid j \in [1,p] \}$$
might not be increasing. Therefore, there are redundancies in Inequalities  
~\ref{basic1}. 
Let $j_1$ be the greatest index such that
$$\sum_{i=1}^{j_1} -\lambda_i= \min\{-\sum_{i=1}^j \lambda_i \mid j \in [1,p] 
\}$$
Then we consider $j \geq j_1$. Let $j_2$ be the greatest number such that 
$$\sum_{i=1}^{j_2} -\lambda_i = \min \{-\sum_{i=1}^{j} \lambda_i \mid j \in 
[j_1,p] \}.$$
If $j_2=j_1$, we stop. Otherwise,
we can continue on and define a sequence
$$j_0=0 < j_1 < j_2 < j_3 < \ldots \leq p$$
with
\begin{equation}~\label{lambdase}
0< \sum_{i=1}^{j_1} -\lambda_i < \sum_{i=1}^{j_1} -\lambda_i < \ldots <
\sum_{i=1}^{p} -\lambda_i.
\end{equation}
Our problem is equivalent to finding  $\{ \eta_{k,i} \}$ such that
$$\eta_{k,i} \in [0,1] \ \ \forall k,i$$
$$\sum_{i=1}^{j_s} (\lambda_i+\sum_{k=1}^n \eta_{k,i}) < 0 \qquad (\forall \ \ 
j_s).$$
Once we determine the sequence
$$j_0=0 < j_1 < j_2 < j_3 < \ldots \leq p,$$
we assign numbers in $[0,1]$ to $\eta_{k,i}$
for $ j_{s-1} < i \leq j_s$ such that
\begin{equation}~\label{al1}
\sum_{i=1}^{j_s} \sum_{k=1}^n \eta_{k,i} < -\sum_{i=1}^{j_s} \lambda_i.
\end{equation}
\begin{thm}~\label{gr1}
Suppose $a \in [1,\infty)^p$. Suppose $\phi(b)b^{-\bold{n+1}}$ on $B^+$
is bounded weakly by $b^{\lambda}$ with $\lambda \prec 0$. 
Then the integral $L(a, \phi)$ converges. Furthermore,
$L(a,\phi)$ is weakly bounded by 
$$a^{-\mu}=\prod_{k=1}^n a_k^{-\mu_k}$$
where $\mu_k= \sum_{i=1}^p \eta_{k,i}$ and for each $j_s>0$, $\{ \eta_{k,i} \in 
[0,1] \}$ satisfy
one of the following
\begin{enumerate}
\item 
\begin{equation}~\label{ar2}
\sum_{i=1}^{j_s} (\lambda_i+\sum_{k=1}^n \eta_{k,i}) =0 ;
\end{equation}
\item
\begin{equation}~\label{ar3}
\sum_{i=1}^{j_s} (\lambda_i+\sum_{k=1}^n \eta_{k,i}) < 0;
\,\,\,\, \mbox{and}
\,\,\,\,  \eta_{k,i}=1  \ \ \forall \ k \in [1,n], i \in [j_{s-1}+1, j_s].  
\end{equation}
\end{enumerate}
\end{thm}
Proof: It suffices to show that for any $0< t <1$, $t \eta_{k,i}$ satisfies the 
conditions in Theorem ~\ref{tw1}. Apparently, we have
$$ t \eta_{k,i} \in [0,1] \qquad (\forall \ i, k)$$
and
$$\sum_{i=1}^{j_s} (\lambda_i+\sum_{k=1}^n \eta_{k,i}) \leq 0$$
From (~\ref{lambdase}), for every $s \geq 1$,
$$\sum_{i=1}^{j_s} (\lambda_i+\sum_{k=1}^n t \eta_{k,i}) \leq (1-t) 
\sum_{i=1}^{j_s}\lambda_i < 0 $$
We have shown that (~\ref{tw11}) holds for $j=j_s$.
For $j_{s-1}+1 \leq j \leq j_{s}$, since $\eta_{k,i} \geq 0$,
\begin{equation}
\begin{split}
 & \sum_{i=1}^{j} \sum_{k=1}^n t \eta_{k,i} \\
 \leq & \sum_{i=1}^{j_{s}} \sum_{k=1}^n t \eta_{k,i} \\
< & -\sum_{i=1}^{j_s} \lambda_i \\
\leq  & -\sum_{i=1}^{j} \lambda_i 
\end{split}
\end{equation}
Thus, (~\ref{tw11}) holds for all $1 \leq j \leq p$. By Theorem ~\ref{tw1},
 $L(a, \phi)$ is bounded by $a^{-t \mu}$ with
$\mu_k= \sum_{i=1}^p \eta_{k,i}$. Hence, $L(a, \phi)$ is weakly bounded by
$a^{-\mu}$. Q.E.D. 
\subsection{$L(p,n)$ and Algorithm for $a \in A^+$ }
Theorem ~\ref{gr1} only assumes $a \in [1, \infty)^n$.
Suppose from now on
$$a \in A^+=\{a_1 \geq a_2 \geq \ldots \geq a_n \geq 1 \}.$$ 
In order to gain a better control over $L(a, \phi)$, we just need to assign 
numbers to $\eta_{1,i}$ to make $\mu_1$ as big as possible, then assign numbers 
to $\eta_{2,i}$ to make $\mu_2$ as big as possible and so on. The only 
requirement is either (~\ref{ar2}) or  (~\ref{ar3}). Our algorithm can be stated 
as follows.
\begin{defn}
Fix $j_s$ and assume that $\{\eta_{k,i} \mid i \leq j_{s-1} \}$ are known.
We assign numbers between $0$ and $1$ to $\eta_{k,i}$ for $j_{s-1} < i \leq 
j_{s}$ in the following way. If (~\ref{ar3}) holds, assign
$\eta_{k,i}=1$ for all $k$ and all $j_{s-1}+1 \leq i \leq j_s$. We are done.
If (~\ref{ar2}) holds,
we choose $\{ \eta_{1,i} \mid  j_{s-1}+1 \leq i \leq j_s \}$ satisfying 
(~\ref{ar2}) and maximizing $\sum_{i=j_{s-1}+1}^{j_s} \eta_{1,i}$. The order of 
assigning numbers to $\{\eta_{1,i} \} $ for $j_{s-1} < i \leq j_{s}$ is not of 
our concern. Update (~\ref{ar2}). 
If (~\ref{ar2}) is trivial, we assign zero to the rest of $\{\eta_{k,i} \mid 
j_{s-1}+1 \leq i \leq j_s \}$ and stop. If not, choose $\{\eta_{2,i} \mid j_{s-
1}+1 \leq i \leq j_s \}$ satisfying (~\ref{ar2}) and maximizing $\sum_{i=j_{s-
1}+1}^{j_s} \eta_{2,i}$. Update ~\ref{ar2} and repeat this process. We do this 
for each $j_s$ until we reach $i=p$. Finally, we compute
$$\mu_k=\sum_{i=1}^p \eta_{k,i} \qquad (1 \leq k \leq n)$$
and obtain a unique $\mu$. Write
$$L(p,n)(\lambda)=-\mu.$$
\end{defn}
The domain of $ L(p,n)$ are apparently $p$-dimensional real vectors such that
$$\lambda \prec 0.$$
The range of $L(p,n)$ are  $n$-dimensional real vectors such that
$$\mu \prec 0.$$
$L(p,n)$, in general, does not produce the precise information for the Langlands 
parameters under theta correspondence; but for a special class of 
representations, $L(p,n)$ will be precise.
Now, Theorem ~\ref{gr1} can be restated as follows.
\begin{thm}~\label{gr2}
Suppose $a \in A^+$. Suppose $\phi(b)b^{\bold{-n+1}}$ on $B^+$
is bounded weakly by $b^{\lambda}$ with $\lambda \prec 0$. 
Then the integral $L(a, \lambda)$ converges. Furthermore,
$L(a,\lambda)$ is weakly bounded by 
$a^{\mu}$ for $\mu= L(p,n)(\lambda)$.
\end{thm}
\subsection{ Examples }
Now let us compute a few examples. Suppose $p \leq n$. \\
\\
{\bf Example 1}: For 
$$\lambda=(-\frac{1}{2}, -\frac{3}{2}, \ldots, -p+\frac{1}{2}),$$
$$L(p,n)(\lambda)=(-p+\frac{1}{2}, -p+1+\frac{1}{2}, \ldots, -\frac{1}{2}, 0, 
\ldots, 0).$$
{\bf Example 2}: For
$$\lambda=(-1, -2, \ldots, -p),$$
$$L(p,n)(\lambda)=(-p, -p+1, \ldots, -1, 0, \ldots, 0).$$
{\bf Example 3}: For
$$\lambda=(-\frac{1}{2}, -\frac{3}{2}, \ldots, -n+\frac{1}{2}),$$
$$L(n,p)(\lambda)=(-n+\frac{1}{2}, -n+\frac{3}{2}, \ldots, -n-\frac{1}{2}+p).$$
{\bf Example 4}: For
$$\lambda=(-1,-2,\ldots, -n),$$
$$L(n,p)(\lambda)=(-n,-n+1, \ldots, -n+p-1).$$

\section{Dual Pair $(O(p,q), Sp_{2n}(\mb R))$ and Estimates on $\theta_s(\pi)$}

Let $O(p,q)$ be the orthogonal group preserving the symmetric form defined by
$$I_{p,q}=\left( \begin{array}{clcr} 0_p & 0 & I_p \\ 0 & I_{q-p} & 0 \\ I_p & 0 
& 0_p
 \end{array} \right)$$ 
and $Sp_{2n}(\mb R)$ be the standard symplectic group. 
We define a symplectic form on $V=M(p+q, 2n)$ by
$$\Omega(v_1, v_2)=Trace(v_1 W v_2^t I_{p,q}) \qquad (\forall \ \ v_1,v_2 \in 
V).$$
Now as a dual pair
 in $Sp(V, \Omega)$, $O(p,q)$ acts by left multiplication and $Sp_{2n}(\mb R)$ 
acts 
by (inverse) right multiplication. We denote both
actions on $M(p+q,2n)$ by $m$. 
\subsection{The dual pair representation $\omega(p,q;2n)$}
Let $x_{i,j}$ be the entries in first $n$ columns of $v \in V$ and $y_{i,j}$ be 
the
entries in the second $n$ columns of $v$. Let 
$$X=\{ v \in V \mid y_{i,j}=0 \}, \qquad Y=\{v \in V \mid x_{i,j}=0\}.$$
Then $X$ and $Y$ are both Lagrangian subspaces of $(V, \Omega)$.
We realize the
Schr\"odinger model of $Mp(V, \Omega)$ on 
$L^2(X)$. Let $\mc P(p,q;2n)$ be the Harish-Chandra module. We call the
admissible representation
$$(\omega(p,q;2n), \mc P(p,q;2n))$$
the dual pair representation. \\
\\
Now let $b=\diag(b_1,b_2,\ldots b_p,1,\ldots,1,b_1^{-1},\ldots,b_p^{-1})$. Let
$$B^+=\{b \mid b_1 \geq b_2 \geq \ldots \geq b_p \geq 1\} \subseteq O(p,q) .$$
Let 
$a=\diag(a_1^{-1},a_2^{-1},\ldots,a_n^{-1},a_1,\ldots, a_n)$. Let
$$A^+=\{a \mid a_1 \geq a_2 \geq \ldots \geq a_n \geq 1\}  \subseteq Sp_{2n}(\mb 
R).$$
 For $1\leq j \leq n$, let
\[ m(b)e_{i,j}= \left\{ \begin{array}{ll}
                b_i e_{i,j} & {i=1,\ldots,p} \\
                e_{i,j} & {i=p+1,\ldots,q} \\
                b_i^{-1} e_{i,j} & {i=q+1, \ldots, p+q}
\end{array} \right.
\]
$$ m(a) e_{i,j}=a_j e_{i,j} \qquad (i=1,\ldots,p+q)$$
These formulae indicate that the embedding $m$ of $A$ and $B$ into $GL(X)$ are
simply the left multiplication and the (inverse) right multiplication. In fact,
$$m(ab) e_{i,j}=\left\{ \begin{array}{ll}
                b_i a_j e_{i,j} & {i=1,\ldots,p} \\
                a_j e_{i,j} & {i=p+1,\ldots,q} \\
               b_i^{-1} a_j e_{i,j} & {i=q+1, \ldots, p+q}
\end{array} \right. $$
Let $b(g_1)$ be the middle term of $KAK$ decomposition of $g_1$ with $b(g_1) \in 
B^+$. Let $a(g_2)$ be the middle term of $KAK$ decomposition of $g_2$ with 
$a(g_2) \in A^+$. Observe that 
$$(b_i a_j+b_i^{-1} a_j^{-1})(b_i^{-1} a_j+b_i a_j^{-1})=
(b_i^2+b_i^{-2}
+a_j^2+a_j^{-2}).$$
From Theorem ~\ref{metaplectic}, we obtain
\begin{thm}~\label{h} For any $\phi, \psi \in \mc P(p,q;2n)$, 
$$| (\omega(p,q;2n)(m(ab)) \phi, \psi) | \leq C \prod_{i=1}^p \prod_{j=1}^n 
(b_i^2+b_i^{-2}
+a_j^2+a_j^{-2})^{-\frac{1}{2}} \prod_{j=1}^n (a_j + a_j^{-1})^{-\frac{q-
p}{2}}$$
Furthermore, this estimate holds for $m(g_1 g_2)$ by substituting $b(g_1)$ and 
$a(g_2)$ into
the right hand side.
\end{thm}
We denote
$$\prod_{i=1}^p \prod_{j=1}^n (b_i^2+b_i^{-2}
+a_j^2+a_j^{-2})^{-\frac{1}{2}}$$
by $H(a,b)$.
\subsection{Growth Control on $\theta_s(p,q;2n)(\pi)$}
Let $(\pi, V)$ be an irreducible Harish-Chandra module in $\mc R_{s}(p,q;2n)$.
We are interested in the following integral
$$\int_{MO(p,q)}(\omega(p,q;2n)(g_1g_2) \phi, \psi)(v, \pi(g_1)u) d g_1 \qquad 
(u, v \in V; \psi, \psi \in \mc P(p,q;2n)).$$
Our goal is to control the growth of this integral as a function on 
$MSp_{2n}(\mb R)$.
From Theorem ~\ref{h} and Theorem ~\ref{equ}, we may as well consider 
\begin{equation}~\label{twisted0}
\int_{B^+} \prod_{j=1}^n (a_j + a_j^{-1})^{-\frac{q-p}{2}} H(a,b) b^{\lambda} 
b^{2 \rho_1} \prod_{i=1}^p\frac{d b_i}{ b_i}
\end{equation}
Here $\rho_1$ is the
half sum of the restricted positive roots of
$O(p,q)$:
$$\rho_1=(\frac{p+q-2}{2}, \frac{p+q-4}{2}, \ldots, \frac{q-p}{2})$$
and $(\pi(g_1)u, v)$ is  bounded by a multiple of
$b(g_1)^{\lambda}$. We observe that
$$\prod_{j=1}^n (a_j + a_j^{-1})^{-\frac{\bold q-\bold p}{2}}\int_{B^+} H(a,b) 
b^{\lambda} b^{2 \rho_1} \prod_{i=1}^p \frac{d b_i}{ b_i}
\leq C a(g_2)^{\bold{-\frac{q-p}{2}}} L(a, \lambda+2 \rho_1- \bold 1)$$
From Theorem ~\ref{gr2}, we obtain
\begin{lem} Let $\pi \in \mc R_s(p,q;2n)$. 
Suppose $K$-finite matrix coefficients of $\pi$ are bounded by some $ C 
b(g_1)^{\lambda}$ with
$$\lambda+ 2\rho(O(p,q))-\bold n \prec 0.$$
Then the matrix coefficients of $\theta_s(p,q;2n)(\pi)$ are weakly bounded by
$$a(g_2)^{L(p,n)(\lambda+2 \rho(O(p,q))-\bold n)-\bold{\frac{q-p}{2}}}.$$
\end{lem}
Recall that $\pi \in \mc R_{ss}(p,q;2n)$ if and only if 
$$\Re(v)-(\bold{n}-\frac{\bold{p+q}}{2})+\rho(O(p,q)) \preceq 0  $$
for every leading exponent $v$ of $\pi$. Take
$$\lambda=\bold{n}-\frac{\bold{p+q}}{2}-\rho(O(p,q))+(\delta,0,\ldots, 0)$$
with $\delta$ a small positive number.
Then matrix coefficients of $\pi$ are bounded by multiples of 
$b(g_1)^{\lambda}$.
\begin{equation}
\begin{split}
 & L(p,n)(\lambda+2 \rho(O(p,q))-\bold{n}) \\
= & L(p,n)(-\frac{\bold{p+q}}{2}+\rho(O(p,q))+(\delta,0, \ldots,0)) \\
= & L(p,n)(-1+\delta, -2, \ldots, -p ) \\
= & \{ \begin{array}{clcr} (-p+\delta, -p+1, \ldots, -1, 0, \ldots, 0)&  n \geq p 
\\
(-p+\delta,-p+1, \ldots, -p+n-1) & n < p
\end{array}.
\end{split}
\end{equation}
From Lemma ~\ref{delta}, we obtain the following theorem
\begin{thm}~\label{thetaspq}
Suppose that $\pi \in \mc R_{ss}(p,q;2n)$. Then
the matrix coefficients of $\theta_s(p,q;2n)(\pi)$ are weakly bounded by
$$a(g_2)^{(-\frac{p+q}{2}, -\frac{p+q-2}{2}, \ldots, -\frac{q-p}{2}, \ldots, -
\frac{q-p}{2})} \qquad  (if \  n \geq p)$$
$$a(g_2)^{(-\frac{p+q}{2},-\frac{p+q-2}{2}, \ldots, -\frac{p+q-2n+2}{2})} \qquad 
(if \  n < p).$$
\end{thm}

\subsection{Growth Control on $\theta(2n;p,q)_s(\pi)$}
Let $(\pi, V)$ be an irreducible Harish-Chandra module in $\mc R_{s}(2n;p,q)$.
We are interested in the following integral
$$\int_{MSp_{2n}(\mb R)}(\omega(p,q;2n)(g_1g_2) \phi, \psi)(v, \pi(g_2)u) d g_2 
\qquad (u, v \in V;\phi, \psi \in \omega(p,q;2n)).$$
Our goal is to control the growth of this integral as a function on $MO(p,q)$.
From Theorem ~\ref{h} and Theorem 8.47 in ~\cite{knapp}, it suffices to consider 
\begin{equation}~\label{twisted0}
\int_{A^+} H(a,b) a^{\lambda} a^{2 \rho_2} \prod_{j=1}^n (a_j + a_j^{-1})^{-
\frac{q-p}{2}} \frac{d a_j}{ a_j}
\end{equation}
Here $\rho_2$ is
the half sum of the restricted positive roots of $Sp_{2n}(\mb R)$:
$$\rho_2=(n, n-1, \ldots, 1)$$
and $(\pi(g_2)u, v)$ is  bounded by a multiple of
$a(g_2)^{\lambda}$. Apparently, the integral ~\ref{twisted0} can be controlled 
by $C L(a, \lambda-\bold{\frac{q-p}{2}}-\bold{1}+ 2 \rho_2)$. From Theorem 
~\ref{gr2}, we obtain
\begin{lem}~\label{control}
Suppose that $\pi \in \mc R_s(2n;p,q)$, i.e., the matrix coefficients of $\pi$ 
are  bounded by multiples of $a(g_2)^{\lambda}$ for some
$$\lambda+ 2 \rho_2- \frac{\bold{p+q}}{2} \prec 0.$$ Then
the matrix coefficients of $\theta_s(2n;p,q)(\pi)$ are weakly bounded by
$$b(g_1)^{L(n,p)(\lambda+2 \rho_2-\frac{\bold{p+q}}{2})}.$$
\end{lem}
Recall that the representation $\pi$ is in $\mc R_{ss}(2n;p,q)$ if and only if
$$\Re(v)+\bold{n+1}+\rho_2-\frac{\bold{p+q}}{2} \preceq 0$$
for every leading exponent $v$ of $\pi$. Now let
$$\lambda=\bold{-n-1}-\rho_2+\frac{\bold{p+q}}{2}+(\delta,0, \ldots, 0)$$
where $\delta$ is a small positive number. Then the matrix coefficients of
$\pi$ are bounded by multiples of $a(g_2)^{\lambda}$ and
 $$\lambda+2 \rho_2-\frac{\bold{p+q}}{2}=\bold{-n-1}+\rho_2+\delta=(-1+\delta,-
2, \ldots -n).$$
Therefore
$$L(n,p)(\lambda+2 \rho_2-\frac{p+q}{2})=(-n+\delta,-n+1, \ldots, -1,0, \ldots, 
0) \qquad (p>n)$$
$$L(n,p)(\lambda+2 \rho_2-\frac{p+q}{2})=(-n+\delta,-n+1, \ldots, -n+p-1) \qquad 
(p \leq n)$$
From Lemma ~\ref{delta}, we obtain
\begin{thm}~\label{matsp}  Suppose that $\pi$ is in $\mc R_{ss}(2n;p,q)$.
Then matrix coefficients of $\theta(2n;p,q)_s(\pi)$ is weakly bounded by
$$b(g_1)^{(-n,-n+1,\ldots,-1,0,\ldots,0)} \qquad (p>n)$$
$$b(g_1)^{(-n,-n+1, \ldots, -n+p-1)} \qquad (p \leq n).$$
\end{thm}

\subsection{Applications to Unitary Representations}
We may now combine our results from ~\cite{unit} with the results we established 
in the previous two sections.
Let us start with a unitary representation in $\mc R_{ss}(p,q;2n)$. 
\begin{thm} Suppose $p+q \leq 2n+1$. Suppose $\pi$ is a unitary representation 
in $\mc R_{ss}(p,q;2n)$ and $(,)_{\pi}$ is nonvanishing. Then 
$\theta_s(p,q;2n)(\pi)$ is unitary. Furthermore, the matrix coefficients of
$\theta(p,q;2n)(\pi)$ is weakly bounded by
$$a(g_2)^{(\overbrace{-\frac{p+q}{2}, -\frac{p+q-2}{2}, \ldots, -\frac{q-p}{2}-
1}^p,\overbrace{-\frac{q-p}{2}, \ldots, -\frac{q-p}{2}}^{n-p})}$$
\end{thm}
In ~\cite{unit}, we have proved that for $p+q$ odd we can loose our restrictions 
from $\mc R_{ss}(p,q;2n)$ a little bit and unitarity still holds for 
$\theta_s(p,q;2n)(\pi)$.
The precise statement can be stated as follows.
\begin{thm}~\label{odd} Suppose $p+q \leq 2n+1$ and $p+q$ is odd. Suppose $\pi$ 
is a unitary representation in $\mc R_{s}(p,q;2n)$ such that each leading 
exponent $v$ of $\pi$ satisfies
$$\Re(v)-(\bold{n-\frac{p+q-1}{2}})+\rho(O(p,q)) \preceq 0  .$$
If $(,)_{\pi}$ is nonvanishing, then $\theta_s(p,q;2n)(\pi)$ is unitary. 
Furthermore, the matrix coefficients of
$\theta_s(p,q;2n)(\pi)$ is weakly bounded by
$$a(g_2)^{(\overbrace{-\frac{p+q-1}{2}, -\frac{p+q-3}{2}, \ldots, -\frac{q-
p+1}{2}}^p, \overbrace{-\frac{q-p}{2} \ldots, -\frac{q-p}{2}}^{n-p})}$$
\end{thm}
Similarly, we obtain the following theorem regarding $\theta_s(2n;p,q)(\pi)$.
\begin{thm}
Suppose that $n < p \leq q$. Suppose that $\pi$ is a unitary representation in
$\mc R_{ss}(2n;p,q)$.
If $(,)_{\pi}$ is nonvanishing,
then $\theta_s(2n; p,q)(\pi)$ is unitary. Furthermore,
the matrix coefficients of $\theta_s(2n;p,q)(\pi)$ are weakly bounded by
$$b(g_1)^{(\overbrace{-n,-n+1,\ldots,-1}^n,\overbrace{0,\ldots,0}^{p-n})}.$$
\end{thm}

\section{The Idea of Quantum Induction}
In this section, we will define quantum induction first. Then we compute the 
infinitesimal characters of quantum induced representations. Finally, we give 
some indication how the limit of quantum induction will become parabolic 
induction.
\subsection{Quantum Induction on Orthogonal Group}
 Consider the composition of $\theta_s(p,q;2n)$ with 
$\theta_s(2n;p^{\p},q^{\p})$.
Suppose $\pi \in \mc R_{ss}(p,q;2n)$ and $p+q \leq 2n+1$.  If $(,)_{\pi}$ is 
nonvanishing, then $\theta_s(p,q;2n)(\pi)$ is unitary and its leading exponents 
satisfy
$$ \Re(v) \preceq (\overbrace{-\frac{p+q}{2},-\frac{p+q-2}{2}, \ldots, -\frac{q-
p+2}{2}}^p,\overbrace{-\frac{q-p}{2}, \ldots -\frac{q-p}{2}}^{n-p}).$$
The representation $\theta_s(p,q;2n)(\pi)$ is in $\mc R_{ss}(2n;p^{\p}, q^{\p})$ 
if 
$$(\overbrace{-\frac{p+q}{2},-\frac{p+q-2}{2}, \ldots, -\frac{q-p+2}{2}}^p, 
\overbrace{-\frac{q-p}{2}, \ldots -\frac{q-p}{2}}^{n-
p})+\bold{(n+1)}+\rho(Sp_{2n}(\mb R))-\frac{\bold{p^{\p}+q^{\p}}}{2} \preceq 0. 
$$
This is true if and only if
$$-\frac{p+q}{2}+n+1+n-\frac{p^{\p}+q^{\p}}{2} \leq 0.$$
We obtain
\begin{thm}~\label{quan1} Suppose 
$$q^{\p} \geq p^{\p} > n$$
$$p^{\p}+q^{\p}-2n \geq 2n-(p+q)+2 \geq 1$$
$$p+q =p^{\p}+q^{\p} \qquad \pmod 2.$$
 Let $\pi$ be an irreducible unitary representation in $\mc R_{ss}(p,q;2n)$.
Suppose that $(,)_{\pi}$ does not vanish. Then $\theta_s(p,q;2n)(\pi)$ is 
unitary
and $$\theta_s(p,q;2n)(\pi) \in \mc R_{ss}(2n;p^{\p},q^{\p}).$$
Furthermore, $\theta_s(2n; p^{\p},q^{\p})\theta_s(p,q;2n)(\pi)$
is either a unitary representation or the NULL representation.
\end{thm}
\begin{defn} 
Let $\pi$ be a unitary representation in $\mc R_{ss}(p,q;2n)$. Suppose that
$$q^{\p} \geq p^{\p} > n$$
$$p^{\p}+q^{\p}-2n \geq 2n-(p+q)+2 \geq 1$$
$$p+q =p^{\p}+q^{\p} \qquad (\mod 2)$$
We call
$$Q(p,q;2n;p^{\p},q^{\p}): \pi \rightarrow \theta_s(2n; 
p^{\p},q^{\p})\theta_s(p,q;2n)(\pi)$$
the (one-step) quantum induction.
\end{defn}
If one of $(,)_{\pi}$ and $(,)_{\theta(p,q;2n)(\pi)}$ vanishes, we define our 
quantum induction $Q(p,q;2n;p^{\p}, q^{\p})(\pi)$ to be the NULL representation. 
\subsection{Quantum Induction on Symplectic Group}
Next, we consider the composition of $\theta_s(2n;p,q)$ with 
$\theta_s(p,q;2n^{\p})$.
Suppose $n < p \leq q$.
Let $\pi$ be a unitary representation in $\mc R_{ss}(p,q;2n)$. Suppose 
$(,)_{\pi}$ is not vanishing. Then the leading exponents of $\theta(2n;p,q)$ 
satisfy
$$\Re(v) \preceq (-n,-n+1,\ldots,-1,0,\ldots,0).$$
Therefore, $\theta(2n;p,q)$ is in $\mc R_{ss}(MO(p,q), \omega(p,q;2n^{\p}))$ if 
$$(-n,-n+1,\ldots,-1,0,\ldots,0) -\bold{n^{\p}}+ 
\frac{\bold{p+q}}{2}+\rho(O(p,q)) \preceq 0$$
This is true if and only if 
$$-n-n^{\p}+p+q-1 \leq 0.$$
\begin{thm}~\label{quan2}
Suppose $2n^{\p}-p-q \geq p+q-2n-2$ and $n < p \leq q$.
Suppose $\pi$ is a unitary representation in $\mc R_{ss}(2n;p,q)$. If 
$(,)_{\pi}$ does not vanish, then $\theta_s(2n;p,q)(\pi)$ is unitary and it is 
in $\mc R_{ss}(p,q;2n^{\p})$.
Furthermore, $\theta_s(p,q;2n^{\p})\theta_s(2n;p,q)(\pi)$ is a unitary 
representation or the NULL representation.
\end{thm}
\begin{defn}
Let $p,q,n,n^{\p}$ be nonnegative integers such that
$$n < p \leq q$$
$$p+q-2n-2 \leq 2n^{\p}-p-q$$
Let $\pi$ be a unitary representation in $\mc R_{ss}(2n;p,q)$. We call
$$Q(2n;p,q;2n^{\p}): \pi \rightarrow 
\theta_s(p,q;2n^{\p})\theta_s(2n;p,q)(\pi)$$
the (one-step) quantum induction.
\end{defn}
If one of $(,)_{\pi}$ and $(,)_{\theta_s(2n;p,q)(\pi)}$ vanishes, we define our 
quantum induction $Q(2n; p,q;2n^{\p})(\pi)$ to be $0$. Thus the domain of our 
quantum induction is $\mc R_{ss}(2n;p,q)$. 
\subsection{Quantum Inductions}
We can further define $2$-step quantum induction and so on. The general quantum 
induction
$$Q(p_1,q_1; 2n_1; p_2,q_2; 2n_2; \ldots)(\pi)$$
is defined as the composition of $\theta_s$,
under the following
conditions:
\begin{enumerate}
\item Initial Conditions: \\
$p_1+q_1 \leq 2n_1+1$. \\
$\pi$ is a unitary representation in $\mc R_{ss}(p_1,q_2; 2n_1)$, i.e., its 
leading exponents satisfy
$$\Re(v)-\bold{n_1}+\frac{\bold{p_1+q_1}}{2}+\rho(O(p_1,q_1)) \preceq 0$$
 \item Inductive Conditions: $\forall \ \ j$, \\
$$ n_j < p_{j+1} \leq q_{j+1}$$
$$ p_{j+1}+q_{j+1}-2 n_{j} \leq 2n_{j+1}-p_{j+1}-q_{j+1}+2$$
$$ 2 n_j -p_j-q_j+2 \leq p_{j+1}+q_{j+1} - 2 n_j$$
$$ p_j+q_j \equiv p_{j+1}+q_{j+1} \qquad (\mod 2).$$
\end{enumerate}
\begin{thm} The representation
$$Q(p_1,q_1; 2n_1; p_2,q_2; 2n_2; \ldots)(\pi)$$
is either an irreducible unitary representation or the NULL representation.
\end{thm}

The general quantum induction
$$Q(2n_1; p_1,q_1; 2n_2; p_2,q_2; 2n_3; \ldots)(\pi)$$
is defined as the composition of $\theta_s$ under the following
conditions:
\begin{enumerate}
\item Initial Conditions: \\
$n_1 < p_1 \leq q_1$ \\
$\pi$ is a unitary representation in $\mc R_{ss}(2n_1; p_1,q_1)$, i.e., its 
leading exponents satisfy
$$\Re(v)-\frac{\bold{p_1+q_1}}{2}+\bold{n+1}+\rho(Sp_{2n_1}(\mb R)) \preceq 0$$
 \item Inductive Conditions: $\forall \ \ j$, \\
$$ n_j < p_{j} \leq q_{j}$$
$$ p_j+q_j-2 n_{j} \leq 2n_{j+1}-p_j-q_j+2$$
$$ 2 n_{j+1} -p_j-q_j+2 \leq p_{j+1}+q_{j+1} - 2 n_{j+1}$$
$$ p_j+q_j \equiv p_{j+1}+q_{j+1} \qquad (\mod 2).$$
\end{enumerate}
\begin{thm}The representation
$$Q(2n_1; p_1,q_1; 2n_2; p_2,q_2; 2n_3; \ldots)(\pi)$$
is either an irreducible unitary representation or the NULL representation.
\end{thm}

Our inductive conditions are natural within the frame work of orbit method (see 
~\cite{vogan3}, ~\cite{thesis}, ~\cite{pan}, ~\cite{pr1}). The nonvanishing of 
$\theta_s$ has been
studied in ~\cite{thesis} and ~\cite{non1}. It can be assumed as a working 
hypothesis in the framework of quantum induction. Notice that $Q$ is defined 
as a composition of $\theta_s$. Thus, it is not known that $Q$ is exactly the 
composition of theta correspondences over $\mb R$. This problem hinges on one
earlier problem mentioned by Jian-Shu Li (see \cite{li1}): \\
\\
Is $(,)_{\pi}$ nonvanishing if $\pi \in \mc R(MG_1,MG_2) \cap \mc R_{s}(MG_1, 
MG_2)$? \\
\\
Our result in ~\cite{theta} which is derived from Howe's results in ~\cite{howe} 
confirms the converse:\\
\\
$\pi$ is in $\mc R(MG_1,MG_2)$ if $(,)_{\pi}$ does not vanish.  \\
\\
Therefore, if $Q(*)(\pi) \neq 0$, $Q(*)$ is the composition of $\theta$. 
\subsection{Infinitesimal Characters}
Infinitesimal characters under theta correspondence were studied by Przebinda 
(\cite{pr}). We denote the infinitesimal character of an irreducible 
representation $\pi$ by $\mc I(\pi)$.
Przebinda's result can be stated as follows.
\begin{thm}[Przebinda] 
\begin{enumerate}
\item Suppose $p+q < 2n+1$. Then 
$$\mc I(\theta(p,q;2n)(\pi))=\mc I(\pi) \oplus (n-\frac{p+q}{2}, n-
\frac{p+q}{2}-1, \ldots, 1+[\frac{p+q}{2}]-\frac{p+q}{2}).$$
\item Suppose $2n+1 < p+q$. Then $$\mc I(\theta(2n;p,q)(\pi)) =\mc I(\pi) \oplus 
(\frac{p+q}{2}-n-1, \frac{p+q}{2}-n-2,\ldots, \frac{p+q}{2}-[\frac{p+q}{2}]).$$
\item Suppose $p+q=2n$ or $p+q=2n+1$. Then $\mc I(\theta(p,q;2n)(\pi))=\mc 
I(\pi)$.
\end{enumerate}
\end{thm}
Now we can compute the infinitesimal character under quantum induction.
\begin{cor} Suppose $Q(*)(\pi) \neq 0$.
\begin{enumerate}
\item If $p+q$ is even, then 
$$\mc I(Q(2n;p,q;2n^{\p})(\pi))=\mc I(\pi) \oplus (\frac{p+q}{2}-n-1, 
\frac{p+q}{2}-n-2, \ldots, 0) \oplus (n^{\p}-\frac{p+q}{2}, n^{\p}-
\frac{p+q}{2}-1, \ldots, 1).$$
\item If $p+q$ is odd, then 
$$\mc I(Q(2n;p,q;2n^{\p})(\pi))=\mc I(\pi) \oplus
(\frac{p+q}{2}-n-1, \frac{p+q}{2}-n-2, \ldots, \frac{1}{2}) \oplus (n^{\p}-
\frac{p+q}{2}, n^{\p}-\frac{p+q}{2}-1, \ldots, \frac{1}{2}).$$
\item If $p+q$ is even, then 
$$\mc I(Q(p,q;2n;p^{\p},q^{\p})(\pi))=\mc I(\pi) \oplus (n-\frac{p+q}{2}, n-
\frac{p+q}{2}-1, \ldots, 1) \oplus (\frac{p^{\p}+q^{\p}}{2}-n-1, 
\frac{p^{\p}+q^{\p}}{2}-n-2, \ldots, 0).$$
\item If $p+q$ is odd, then 
$$\mc I(Q(p,q;2n;p^{\p},q^{\p})(\pi))=\mc I(\pi) \oplus (n-\frac{p+q}{2}, n-
\frac{p+q}{2}-1, \ldots, \frac{1}{2}) \oplus (\frac{p^{\p}+q^{\p}}{2}-n-1, 
\frac{p^{\p}+q^{\p}}{2}-n-2, \ldots, \frac{1}{2}).$$
\end{enumerate}
\end{cor}
We shall now take a look at some "limit" cases under quantum induction. \\
\\
{\bf Example I}: $p+q+p^{\p}+q^{\p}=4n+2$. \\
\\
In this case,
$$n-\frac{p+q}{2}=\frac{p^{\p}+q^{\p}}{2}-n-1.$$
Therefore,
$$\mc I(Q(p,q;2n;p^{\p},q^{\p})(\pi))=\mc I(\pi) \oplus \overbrace{ (n-
\frac{p+q}{2},
n-\frac{p+q}{2}-1, \ldots, 1+\frac{p+q}{2}-n, \frac{p+q}{2}-n)}^{2n-p-q+1}.$$
{\bf Example II}: $2n-p-q+2=p^{\p}+q^{\p}-2n$ and $p-p^{\p}=q-q^{\p}$.\\
\\
Notice first that
$$p^{\p}-p+q^{\p}-q=(p^{\p}+q^{\p})-(p+q)=4n+2-2(p+q).$$
Therefore
$$\frac{p^{\p}-p}{2}=\frac{p^{\p}-p+q^{\p}+q}{4}=n-\frac{p+q}{2}+\frac{1}{2}$$
Recall from Prop 8.22  ~\cite{knapp} 
\begin{equation}
\begin{split}
  & \mc I(Ind_{SO_0(p,q) GL_0(p^{\p}-p) N}^{SO_0(p^{\p},q^{\p})}(\pi \otimes 1)) 
\\ = &
\mc I(\pi \otimes 1) \\
= & \mc I(\pi) \oplus (\frac{p^{\p}-p-1}{2}, \frac{p^{\p}-p-3}{2}, \ldots, -
\frac{p^{\p}-p-3}{2}, -\frac{p^{\p}-p-1}{2}) \\
= & \mc I(\pi) \oplus (n-\frac{p+q}{2},n-\frac{p+q}{2}-1, \ldots, 
1+\frac{p+q}{2}-n, \frac{p+q}{2}-n) \\
= & \mc I(Q(p,q;2n;p^{\p},q^{\p})(\pi)).
\end{split}
\end{equation}
This suggests that $Q(p,q;2n;p^{\p},q^{\p})(\pi)$ as a representation of 
$SO_0(p,q)$ can be decomposed as direct sum of some parabolically induced 
unitary representation (see Conjecture I). \\
\\
{\bf Example III}: $n+n^{\p}+1=p+q$. \\
\\
In this case,
$$\frac{p+q}{2}-n-1=n^{\p}-\frac{p+q}{2},$$
$$\frac{n^{\p}-n-1}{2}=\frac{p+q}{2}-n-1.$$
From Prop. 8.22 (~\cite{knapp}) and the Corollary,
\begin{equation}
\begin{split}
 &\mc I(Ind_{Sp_{2n}(\mb R) GL(n^{\p}-n)N}^{Sp_{2n^{\p}}(\mb R)} (\pi \otimes 
1)) \\
= & \mc I(\pi) \oplus (\frac{n^{\p}-n-1}{2}, \frac{n^{\p}-n-3}{2}, \ldots,
-\frac{n^{\p}-n-3}{2}, -\frac{n^{\p}-n-1}{2}) \\
= & \mc I(\pi) \oplus
(\frac{p+q}{2}-n-1, \frac{p+q}{2}-n-2, \ldots, -\frac{p+q}{2}+n+2, -
\frac{p+q}{2}+n+1) \\
= & \mc I(Q(2n;p^{\p},q^{\p};2n^{\p})(\pi))
\end{split}
\end{equation}
This suggests that $Q(2n;p,q;2n^{\p})(\pi)$ can be obtained as subfactors of 
certain
parabolic induced representation. We prove this connection in ~\cite{quan}.\\
\\
Let me make some final remarks concerning the definition of quantum induction
$Q$. Notice that $ Q(p,q;2n;p^{\p},q^{\p})(\pi)$ contains distributions of the 
following form
\begin{equation}
\begin{split}
 & \int_{MSp_{2n}(\mb R)} \omega(p^{\p},q^{\p};2n)(g_1) \phi_1 \otimes
\int_{MO(p,q)} \omega(p,q;2n)^c (g_1g_2) \phi_2 \otimes \pi(g_2) v d g_2 dg_1 \\
= & \int_{MO(p,q)} \omega(p,q;2n)^c(g_2) [\int_{MSp_{2n}(\mb R)} 
\omega(p^{\p}+q, q^{\p}+p;2n)(g_1) (\phi_1 \otimes \phi_2 ) d g_1] \otimes 
\pi(g_2) v.
\end{split}
\end{equation}
Our discussions in this paper guaranteed absolute integrability of this 
integral. Notice that the vectors in $[*]$ are in
$\theta(2n;p^{\p}+q,q+p^{\p})(1)$.

\begin{defn} Suppose $p^{\p}+q \geq 2n$, $q^{\p}+p \geq 2n$ and
$p+q+p^{\p}+q^{\p}$ is even.
Consider the dual pair $(O(p^{\p}+q,q^{\p}+p), Sp_{2n}(\mb R))$. This is a dual 
pair in the stable range (~\cite{li2}, ~\cite{howesmall}).
Then $\theta(2n;p^{\p}+q,q^{\p}+p)(1)$ is an unitary representation of 
$MO(p^{\p}+q,q^{\p}+p)$ (see ~\cite{li2}, ~\cite{hz}). Let $O(p,q)$ and 
$O(p^{\p},q^{\p})$ be embedded diagonally in $O(p^{\p}+q,q^{\p}+p)$. Let $\pi 
\in \Pi(MO(p,q))$. Formally define a Hermitian form $(,)$ on 
$\theta(2n;p^{\p}+q,q^{\p}+p)(1) \otimes \pi$ by integrating the matrix 
coefficients of $\theta(2n;p^{\p}+q,q^{\p}+p)(1)$ against
the matrix coefficients of $\pi$ over $MO(p,q)$ as in (~\ref{avera}). 
Suppose that $(,)$ converges. Define $\mc Q(p,q;2n;p^{\p}, q^{\p})(\pi)$ to be 
$\theta(2n;p^{\p}+q,q^{\p}+p)(1) \otimes \pi$ modulo the radical of $(,)$.  $\mc 
Q(p,q;2n;p^{\p},q^{\p})(\pi)$ is thus a representation of $MO(p^{\p},q^{\p})$.
\end{defn}
One must assume that $p^{\p}+q^{\p} \equiv p+q \pmod 2$. Otherwise, 
$\theta(2n;p^{\p}+q, q^{\p}+p)(1)=0$.
$\mc Q$ can be regarded as a more general definition of quantum induction. It is 
no longer clear that $\mc Q$ preserves unitarity.
\begin{thm} Under the assumptions from Theorem ~\ref{quan1}, 
$$\mc Q(p,q;2n;p^{\p},q^{\p})(\pi) \cong Q(p,q;2n;p^{\p},q^{\p})(\pi).$$
\end{thm}

 Similarly, one can define nonunitary quantum 
induction $\mc Q(2n;p,q;2n^{\p})(\pi)$. 
\begin{defn}
Suppose that $p+q \leq n+n^{\p}+1$. Consider the dual pair $(O(p,q), 
Sp_{2n+2n^{\p}}(\mb R))$. Then $\theta(p,q;2n+2n^{\p})(1)$ is a unitary 
representation of $MSp_{2n+2n^{\p}}(\mb R)$ (see ~\cite{howesmall}, ~\cite{li2}, 
~\cite{pr2}). Let $\pi \in \Pi(MSp_{2n}(\mb R))$. Formally define a Hermitian 
form $(,)$ on $\theta(p,q;2n+2n^{\p})(1) \otimes \pi$
by integrating the matrix coefficients of $\theta(p,q;2n+2n^{\p})(1)$ against 
the matrix coefficients of $\pi$ as in (~\ref{avera}). Suppose that $(,)$ 
converges. Define $\mc Q(2n;p,q;2n^{\p})(\pi)$ to be $\theta(p,q;2n+2n^{\p})(1) 
\otimes \pi$ modulo the radical of $(,)$. $\mc Q(2n;p,q;2n^{\p})(\pi)$ is a 
representation of $MSp_{2n}(\mb R)$.
\end{defn}
For $p+q$ odd, the $MSp$ in this definition are metaplectic groups. For $p+q$ 
even, the $MSp$ in this definition split (see Lemma ~\ref{msp2n}).
\begin{thm}
Under the assumptions from Theorem ~\ref{quan2}, 
$$\mc Q(2n; p,q;2n^{\p})(\pi) \cong Q(2n;p,q;2n^{\p})(\pi).$$
\end{thm}
There is a good chance that $\mc Q(*)(\pi)$ will be irreducible. \\
\\
Quantum induction fits well with the general philosophy of induction. On the one 
hand, similar to parabolic induced representation $Ind_{P}^G \tau$ whose vectors 
are in
$$Hom_P(C^{\infty}_c(G), \tau),$$
quantum induced $\mc Q(p,q;2n;p^{\p},q^{\p})(\pi)$ lies in
$$Hom_{\f o(p,q), O(p) \times O(q)}(\theta(2n; p^{\p}+q,q^{\p}+p)(1), \pi).$$
On the other hand, $Ind_P^G \tau$ has a nice geometric description. It consists 
of sections of 
 the vector bundle
$$ G \times_P \tau \rightarrow G/P.$$
In contrast, quantum induction does not possess this kind of classical 
interpretation except for some limit case.

\end{document}